\documentclass[article,noinfoline]{imsart}

\usepackage{amsthm,amsmath,amssymb,natbib}

\usepackage{graphicx}
\usepackage{color}

\def\argmax{\mathop{\rm argmax}}
\startlocaldefs

\newcommand{\bdf}{\boldsymbol}
\endlocaldefs

\begin{document}

\begin{frontmatter}

% "Title of the paper"
\title{ Deconvolution with application to estimation of sampling probabilities
and the Horvitz-Thompson estimator.  }
\runtitle{Deconvolution Horvitz-Thompson}

% indicate corresponding author with \corref{}
% \author{\fnms{John} \snm{Smith}\corref{}\ead[label=e1]{smith@foo.com}\thanksref{t1}}
% \thankstext{t1}{Thanks to somebody}
% \address{line 1\\ line 2\\ printead{e1}}
% \affiliation{Some University}

%\author{\fnms{Lawrence D.} \snm{Brown,}\thanksref{t1}\ead[label=e1]{lbrown@wharton.upenn.edu}}
%\address{University of Pennsylvania; \printead{e1}}
%\thankstext{t1}{Research supported by an NSF grant.}
%\affiliation{University of Pennsylvania}
\author{\fnms{Eitan} \snm{Greenshtein}\ead[label=e2]{eitan.greenshtein@gmail.com}}
\address{Israel Central Bureau of Statistics; \printead{e2}}
\affiliation{Israel Central Bureau of Statistics}
\and
\author{\fnms{Theodor } \snm{Itskov} \ead[label=e3]{itsmatis@gmail.com}}
\address{   Israel Central Bureau of Statistics ; \printead{e3}}
%\thankstext{t3}{Research supported by an ISF grant.}
\affiliation{Israel Central Bureau of Statistics; \printead{e3}}

\runauthor{Greenshtein, Itskov}

\end{frontmatter}

%\begin{document}

%\title{Deconvolution estimation of sampling probabilities,
%and a modification of the Horvitz-Thompson estimator.}

%\maketitle

{\bf Abstract.}
We elaborate on a deconvolution method, used to estimate the empirical distribution of unknown parameters, as suggested recently by Efron (2013). It is applied to estimating the empirical distribution of the
`sampling probabilities' of $m$ sampled items. The estimated empirical distribution is used to modify the Horvitz-Thompson estimator. The performance of the modified Horvitz-Thompson estimator is studied in two
examples. In one example the sampling probabilities are estimated based on the
number of visits until a response was obtained. 
The  other example is based on real data from panel sampling, where
in four consecutive months there are corresponding four attempts
to interview each member in a panel. The sampling probabilities are estimated based on the number of successful attempts. 

We also discuss briefly, further applications of deconvolution, including estimation of False discovery rate.

\section{Introduction and Preliminaries }

In this paper we suggest a deconvolution method for estimating
the empirical distribution of  unknown parameters that correspond to
a given set of independent observations. Our motivation and our examples, involve the estimation of sampling/response probabilities of  a given set of sampled items. Such an estimation leads to a modified Horvitz-Thompson estimator.

Given a population of  $N$ items with corresponding measurements 
$X_1,...,X_N$, it is desired to estimate the total
$T=\sum^N_{i=1} X_i$, based on randomly  sampled 
items ${\cal {\bdf{S}}}=\{i_1,...,i_m \}$.
Denote $I_i$ the indicator of the event item $i$ was sampled. 
Let $p_i=E I_i$, then the Horvitz-Thompson estimator of T is
$$\hat{T}=\sum_{i \in S} \frac{X_i}{p_i}.$$ Writing $\hat{T}$ as $\sum_{i=1}^N (X_i/p_i) I_i$, it is immediate that $E \hat{T}=T$,
when $min_i p_i>0$.

It is often the case that $p_i, \; i=i_1,...,i_m$, are unknown and 
thus the above estimator can not be  derived; a way for
approximating it is developed in this paper.  
We will suggest a deconvolution method that will estimate functionals of
$p_i, \; i=i_1,...,i_m$, based on the `sampling effort'. Our main example, which also serves for illustrating our general 
deconvolution approach,
 is where the sampling effort is related to the number of `visits'/'calls' required to get a response. The other example
 is based on real data that was obtained in the Labor Force Survey
in Israel, this is a panel sampling where each panel is investigated in four consecutive months.
The response probabilities
(or, required `sampling effort') are estimated based on the number
of responses that were obtained in the four months.

A general reference to sampling is, e.g., Lohr(2009). A reference for missing data and non-response issues, that are involved in our  examples is, e.g., Little and Rubin (2002).

Consider a typical scenario where  item $i$, is selected with known probability $\pi_i$ to a certain list $L$ of size $I$, however, once  item  $i$ is in the list, there is a probability $p^*_i \leq 1$ to get a response from that item when approached for an interview. Denote by $R_i$ the event: `a response
was obtained from item $i$'. The probabilistic model we assume is where   
$$p_i= P(I_i=1)=P(R_i \cap (i\in L))=P(R_i|i \in L) P(i \in L)=p^*_i \pi_i.$$
Note, $p_i$ is unknown since  $p_i^*$ is unknown. 

There are various  ways  for estimating the probabilities $p_i$,
all are based on the assumption that items with similar certain characteristics 
(say, gender, age), have equal response probabilities. Of course, such an
assumption, even if often helpful,
can not be  right in reasonable setups, even if we try to estimate
the response probabilities based on a rich set of characteristics.

We will suggest an approach for estimating
(functionals of) the response probabilities $p_i^*$,
based on the `effort' made to get a response.
Suppose our policy is to make $M_0$ attempts in order to obtain a response, however obviously if a response was obtained in the $j<M_0$ attempt, no further attempts are made. Let $M_i$
be the number of attempts made until a response is obtained. We model
$M_i$ as a Geometric random variable $M_i \sim Ge(\tilde{p}_i)$, assume 
$0<\min_i 
\tilde{p}_i$. Note,
for the items in the list we observe only
the corresponding truncated observations denoted $Y_i$, where
$$P(Y_i=j)= \frac{ (1-\tilde{p}_i)^{j-1}\tilde{p}_i  }
{ 1- (1-\tilde{p}_i)^{M_0} }, \; j=1,...M_0.$$

For simplicity we assume that the measurements
$X_i$, may have only the values $0$ or $1$, though the same treatment may be given to any r.v., with finite number of possible outcomes.
Also, we simplify by treating the case $$\pi_i \equiv I/N.$$ 
Let $m$ be the random size of the set ${\cal {\bdf{S}}}$ of
 sampled items (i.e., those in the list which responded), $m_1$ out of those $m$ items have a corresponding value $X_i=1$,
and $m_0$ have a corresponding value $X_i=0$, $m_1+m_0=m$.

The response probabilities $p_i^*$, i.e., the probability of a response
from item $i, i=1,...,N$ conditional that it was selected to the list,
are equal to $$p_i^*= 1-(1-\tilde{p}_i)^{M_0},$$ according to our model.

We re-index the response probabilities $p_i^*$ that correspond
to sampled   items with $X_i=1$, by $p_{1,i}^*, \; i=1,...,m_1$ and 
those that correspond to sampled   items with $X_i=0$, by $p_{0,i}^*, \; i=1,...,m_0$, similarly we write $\pi_{1,i}$, $Y_{1,i}$, etc. 

The unbiased Horvitz-Thompson estimator  (if $p_{1,i}^*$, $i=1,...,m_1$ were known!) is:
\begin{equation}
\label{eqn:HT}
\sum_{i=1}^{m_1} \frac{1}{\pi_{1,i} p_{1,i}^*} =\frac{N}{I} \sum_{i=1}^{m_1} \frac{1}{p_{1,i}^*}
\equiv \frac{N}{I} \theta, \end{equation}
we will later refer to the last estimator as the `oracle-estimator'.

Our goal is to estimate (the random quantity)
$\theta$, by an appropriate empirical-Bayes
or deconvolution estimator
$\hat{\theta}$ that will yield a Modified Horvitz-Thompson estimator 
of the form $$\hat{T}_{MHT}=\frac{N}{I} \hat{\theta}.$$ 

A standard approach, for estimating $\theta$, is to estimate $p_{1,i}^*$, $i=1,...,m_1$, e.g., by an mle $\hat{p}^*_{1,i}=\hat{p}^*_{1,i}(Y_{1,i})$ and plug-in to the above expression for $\theta$. This approach would yield a biased estimator
even asymptotically. Alternatively, our approach stems from the representation: $$ \frac{\theta}{m_1}= E_G \frac{1}{S},$$ for a random variable $S$ with distribution $G$ that equals to the empirical distribution of 
$\{ p_{1,1}^*,...,p_{1,m_1}^* \}$. Thus, if
we find an estimator $\hat{G}$ that approaches weakly to $G$, we expect
$E_{\hat{G}} 1/S \rightarrow E_G 1/S$, when the support of $G$ is bounded
above 0.

The first study, of mle estimation of a mixing distribution $G$ and weak convergence of the mle estimator, was done by Kiefer and
Wolfowitz (1956). Their setup is in the spirit of empirical Bayes,
treating the unknown parameters as independent realizations from an unknown
$G$. We also discuss the formal setup, 
where we want to estimate the empirical distribution of the sampled parameters; this setup is  in the compound decision spirit. 
The distinction is often suppressed in the literature;
see also the next section about the distinction.

The term deconvolution for our estimation procedure,
might be more appropriate
then Empirical Bayes, since we are only interested in estimating a functional of the empirical distribution of $p_{11}^*,..., p_{1,m_1}^*$,
rather than the vector of individual parameters. 
In non-parametric empirical Bayes,
the estimation of the empirical distribution is often used as a first stage
in estimation of the individual parameters, where the estimated empirical distribution is served as a `prior', see e.g.,
Efron (2013). In some cases the estimation of the empirical distribution may be circumvented,
and the optimal non-parametric empirical Bayes procedure is approximated directly. This is the case in the classical Robbins's procedure for the Poisson problem, see e.g., Brown et. al.  (2013); see, e.g., Brown and Greenshtein (2009) for direct approximation in the normal problem.  In the sequel,
we will further motivate the estimator $E_{\hat{G}} 1/S$ from an EB perspective, see Subsection 2.5.

\bigskip

Our deconvolution is a method for deriving Non-Parametric  Maximum Likelihood Estimator (NPMLE) for empirical distribution.
We use the term deconvolution in a wide sense, that includes identifying mixtures, as studied, e.g., by Lindsay (1995), Lindasy and Roeder (1993).
See further discussion and literature in Brown et. al. (2013) for the Poisson case, Koenker and Mizera (2013),  Jiang and Zhang (2009), Lee, et.al (2013), and, of course, the fore mentioned seminal paper of Kiefer and Wolfowitz (1956).
Our quadratic programming approach, rather than the more common EM-algorithm, is in line with the general suggestion and advocation
of Koenker and Mizera (2013) for the usage of convex optimization. 

\bigskip

In Section 2 we describe our deconvolution method through our
illustrating example, we also discuss its merits and 
various theoretical
issues; 
in Section 3 we present simulation results;
finally, in Section 4, we demonstrate our method  through an analysis
of a real data set obtained in the Labor Force Survey in Israel.

\section{Estimation}

\subsection{ Deconvolution.}

We will consider both, a `Compound Decision (CD) approach',
for estimating the empirical distribution,
denoted $G$, of $p^*_{1,i}, \; i=1,...,m_1$, and the
more common `Empirical Bayes (EB) approach' for estimating a distribution $G$ from which $p_{1i}^*$
are assumed to be an independent sample. 

Our perspective
on the `CD approach' involves
reduction by permutation invariance. When searching for a decision rule which is
permutationally invariant with respect to the observations, we may confine ourselves to functions of
a maximal invariant, e.g., the order statistic of the observations. The corresponding parameter space in the reduced problem, may be taken to be
the order statistic of the parameters. Note, estimating the order statistic
of the parameters is equivalent to estimating their empirical distribution.
We will estimate the order statistics of the parameters both, through
a version of the method of moments and through an approximation of the mle, in the reduced by
invariance problem.

\bigskip

Let $Y_{(1)},...,Y_{(m_1)}$ be the order statistic of 
$Y_{1,1},...,Y_{1,m_1}$, let $p_{(1)},...,p_{(m_1)}$ be the order statistic of $p_{1,1},...,p_{1,m_1}$. The mle in the reduced problem
is:  
$$ \argmax_{p_{(1)},...,p_{(m_1)} } \; P_{p_{(1)},...,p_{(m_1)}} 
{\Large (}Y_{(1)}=y_{(1)},...,Y_{(m_1)}=y_{(m_1)} {\Large )}.$$

An exact computation of the above involves summation 
of the relevant probabilities that correspond to all the
permutations of $y_{(1)},...,y_{(m_1)}$. It is computationally intractable, and an approximation method is described bellow.
The same ideas apply for the `EB-approach'.

\bigskip

In the sequel it is convenient to parametrize according to 
$p_{1i}^* \equiv s_i$, 
rather than according to $\tilde{p}_{1,i}$, $i=1,...,m_1$. 
We write , e.g., 
$P_{s_i}(M_i \leq M_0)=s_i$.

Let  $k \leq m_1$ be the number of distinct values
in the support of the empirical distribution 
of $p_{1,1}^*,..., p_{1,m_1}^*$. 
Denote by $\bdf {g}'=(g_1,...,g_k)$ the discrete density of the empirical distribution
of $\{ p_{1,1}^*,...,p_{1,m_1} ^* \}$.  Denote the points in the support
of the distribution by $\bdf{s}=(s_1,...,s_k)$.

The deconvolution or estimation, of the empirical distribution is through
a modification of a simple discretization method, described and demonstrated
recently by Efron (2013)  in the context of Empirical Bayes.

The idea is explained by first assuming that
the empirical distribution is known. For any point $s_i \in \bdf{s}$, and any
$j=1,...,M_0$, denote $$p_{ji}= P_{s_i} (Y=j);$$ let $P=(p_{ji})$ be the corresponding $M_0 \times k$ matrix.
Given $m_1$ observations,
from the distribution defined by the $\bdf{g}$ mixture of the
distributions that correspond to the parameters in
$\bdf{s}$, a random vector of dimension $M_0$, denoted ${\bf C}$ is induced,
where its $1 \leq j \leq M_0$ coordinate  counts of the number of times
the corresponding (truncated) Geometric variable obtained the value $j$. 
Note, the random vector ${\bf C}$ is maximal invariant in the sense of the first paragraph of this section.
Denote
$\bdf{f}'=(f_1,...,f_{M_0})= E ({\bf C}/m_1)$; note, $\sum f_l=1$.

Then: \begin{equation}  \bdf{f}= P \bdf{g}. \end{equation}
  
This motivates us to approximate the  set $\bdf{s}$ by a `dense'
grid of points $\hat{\bdf{s}}=(\hat{s}_1,...,\hat{s}_\kappa)$, compute the corresponding  $M_0 \times \kappa$ matrix ${P}_{ \bdf{\hat{s}} }$,
and find a vector $\hat{\bdf{g}}$ such that 
\begin{equation}  \bdf{\hat{f}} \approx P_{\hat{\bdf{s}} } \hat{\bdf{g}}; \end{equation}

here $\hat{\bdf{f}}= {\bf C}/m_1$, for a given realization.

More formally, we may reasonably hope that the solution of the following 
quadratic programming problem yields a good approximation of the empirical
distribution of $\{ p_{1,1}^*,...,p_{1,m_1} ^* \}$, for a ` dense
enough' grid
$\hat{\bdf{s}}$.

\begin{equation} \label{eqn:moments} \min_{\bdf{g}} ( \bdf{ \hat{f} }-P_{ \hat{\bdf{s}}} \bdf{g} )' ( \bdf{\hat{f}}-P_{\hat{\bdf{s}}}\bdf{g}), \end{equation}

$\;\;\;\;\;\;$  s.t. $0 \leq g_i \leq 1$, $\sum g_i=1$; \newline
here the support is approximated by the dense grid $\hat{\bf{s}}$, and the solution $\hat{ \bf{g} }$
is an approximation of the mixing probabilities. This approach is a variant
of the method of moments. 
The choice of the support points is rather arbitrary,
the procedure seems only mildly sensitive to that choice. 

An alternative and more elaborated
way is to approximate the mle. Note that, the random vector ${\bf C}$ is asymptotically multivariate normal, with mean $P \bdf{g}$ and a covariance matrix $\Sigma$.  Since $\Sigma$ is singular it is convenient
to consider the vector (and sufficient statistic) ${\bf C}^*$
that consists of the first $(M_0-1)$ coordinates of ${\bf C}$,
whose covariance matrix $\Sigma^*$ is nonsingular and its mean is 
$\bdf{f}^*=(f_1,...,f_{M_0-1})$. Denote 
$\hat{\bdf{f}^*}=(\hat{f}^*_1,...,\hat{f}^*_{M_0-1})$, then,
asymptotically, the mle with respect to the parameters defined by mixtures of the parameters in the  grid 
$\hat{\bdf{s}}$, is determined by the solution of: 

\begin{equation}
 \label{eqn:mle} \min_{ \bdf{g}} ( \bdf{ \hat{f}^* }-P_{ \hat{\bdf{s}}} \bdf{g} )' \Sigma^{*-1} ( \bdf{\hat{f}^*}-P_{\hat{\bdf{s}}}\bdf{g}), \end{equation}

$\;\;\;\;\;\;$  s.t. $0 \leq g_i \leq 1$, $\sum g_i=1$.

Note, in the above $\Sigma^*=\Sigma^*_G$, where $G$ is the empirical
distribution of $s_1,...,s_k$. 
Thus, $\Sigma^*$ is the covariance matrix of the corresponding mixture of multinomials. Hence, the above should be done iteratively, where the mixture distribution is estimated at first stage as in (\ref{eqn:moments}),
yields an initial estimator $\hat{G}$ and a corresponding $\Sigma^*_{\hat{G}}$, which in turn by using
({\ref{eqn:mle}) yields a new estimator for $G$.

\bigskip

{\bf Remark 1 (EB versus CD): } Under the EB-approach
where $p_{1,i}, \; i=1,...,m_1$, are i.i.d realizations from 
some distribution $G$, the covariance matrix $\Sigma^*$ in the above equation corresponds to a multinomial distribution. Under the CD-approach of estimating the empirical distribution 
(or order statistic)
of $p_{1,1}^*,...,p_{1,m_1}^*$, $\Sigma^*$ does not correspond to a multinomial distribution, unless $p_{1,1}^*=...=p_{1,m_1}^*.$
Thus, the two mle approaches yield different estimators, 
that correspond to different estimates for $\Sigma^*_G$. The
possible difference is
demonstrated also
in the following simple Example 1.

{\bf Example 1.}

Consider the most simple example where $Y_i \sim Ber(p_i)$,
$i=1,...,m$. Suppose $\sum_{i=1}^m Y_i=m/2$. Under the EB-approach,
assuming that $p_i$, $i=1,...,m$, is an i.i.d sample from some $G$,
the mle for $G$ is not unique and any  $\hat{G}$ satisfying that $E_{\hat{G}} P=0.5$, is an mle. However, under the CD-approach the unique mle for the empirical distribution of $p_i, \; i=1,...,m,$
has probability 0.5 at 0 and 0.5 at 1.

\bigskip

In our setup the singularity of $P$ and the simulation results indicate that
often the minimum, achieved by the quadratic programming procedure, is
virtually zero, in which case there is no difference between the mle (EB or CD versions) and the
method of moments. 
Thus we used the (simpler) method of moments, specifically,
we used (\ref{eqn:moments}).

The numeric work in this paper was done by applying the 
quadratic programming 
function {\it ipop}, from the R-package {\it kernlab}, Karatzoglou, et. al. (2004).

\subsection{ Generalization and the estimation of $T$.}

In the above we neglected the observations that correspond to
indices $i$ for which the corresponding $X_i=0$. The right way is to
consider both  types of observations together and to utilize the knowledge of the  size $I$ of the initial list.  It is helpful to think of our approach as two stages of inflating the sampled
quantities $m_l, l=0,1$. First
we inflate the observed number $m_l$ of items with $X_i=l$ to get an estimate of their number in the list, and then we inflate 
that estimate again
by multiplying it by $N/I$ to get an estimate of their number in the population.

Denote the vectors ${\bf  \hat{f}_1} \equiv { \bf \hat{f} }$,  
${\bf  {g}_1} \equiv { \bf {g}} \equiv (g_{11},g_{12},...g_{1\kappa})$. 
Denote
the analogous objects that correspond to $X_{0,1},...,X_{0,m_0}$
by  ${\bf  \hat{f}_0}$ and   ${\bf  {g}_0}$. 
Similarly, denote $\Sigma^*_1$ and $\Sigma^*_0$ analogously to the above.
Our quadratic 
programming problem, when considering  both type of
observations  and estimating by $\hat{G}_l$, $l=0,1$,
simultaneously for
both types of empirical distributions
$G_l$, $l=0,1$, is:

\begin{equation}
 \label{eqn:mlegen} \min_{ \bdf{g}_1,\bdf{g}_0} {\Large{[}}
 ( \bdf{ \hat{f}_1^* }-P_{ \hat{\bdf{s}}} \bdf{g}_1 )' \Sigma_1^{*-1} ( \bdf{\hat{f}_1^*}-P_{\hat{\bdf{s}}}\bdf{g}_1)\\ 
 +( \bdf{ \hat{f}_0^* }-P_{ \hat{\bdf{s}}} \bdf{g}_0 )' \Sigma_0^{*-1} ( \bdf{\hat{f}^*_0}-P_{\hat{\bdf{s}}}\bdf{g}_0) {\Large{]}}, 
 \end{equation}

s.t. 

\begin{eqnarray}
&& 0 \leq g_{li} \leq 1, \;l=0,1, \; i=1,...,\kappa ,\\ 
&& \sum g_{1i}=1,  \sum g_{0i}=1, \\
&& m_1 \sum \frac{{g}_{1,i}}{s_i}+m_0 \sum \frac{{g}_{0,i}}{s_i}=I
\end{eqnarray}

The last constraint is applied simultaneously to 
both lists through the known number $I$.

The total
number of units in the population,
with corresponding $X_i=l$, $l=0,1$, is estimated by
\begin{equation}
\label{eqn:TL}
\hat{T}_l= \frac{N}{I} m_1 \sum \frac{\hat{g}_{l,i}}{s_i}= \frac{N}{I}
m_1 E_{\hat{G}_l} \frac{1}{S}=\frac{N}{I} \hat{\theta}_l, \; \; l=0,1.
\end{equation}

The above is trivially generalized to the case where $X_i$
may have  more than two possible outcomes.  As before, this generalization has also a method of moments version.

When there are $L$ possible outcomes and it is desired to
estimate the proportion $pr_l$, $l=1,...,L$ in the population
of units with a corresponding $X=l$, the estimator is:

\begin{equation}
\label{eqn:MHT1}
\hat{pr}_l= \frac{ m_l E_{\hat{G}_l} \frac{1}{S}} 
{\sum m_l E_{\hat{G}_l} \frac{1}{S}}. 
\end{equation} 

\subsection { Calibration.}  

There is a room for further splitting the data beyond just
estimating ${G}_l$, the empirical distribution of response probabilities of units with corresponding $X=l$, $l=1,...,L$. Suppose for example
that for each unit there is an explanotory variable (say, gender)
$W$ that may equal, say, 0 or 1. 
It  may  be helpful to estimate the corresponding distributions
$G_{l,w}$, $l=1,...,L$, $w=0,1$,
and proceed as before. One reason is that $G_{lw}$ might be
more homogeneous and hence may have a good presentation with 
a mixture that has a smaller
number of points in its support, see also the discussion in  subsection 2.5.

However, there could be  another important and good reason for it.
Suppose we know that the proportion in the population of items with
corresponding $W=0$ is 0.5. Then we expect that this is nearly their
proportion in the original list $L$.
Denote by ${\bf g}_{lw}=(g_{lw1}, g_{lw2},...)$, $w=0,1$, $l=1,...,L$ the corresponding
densities, and by $\hat{g}_{lwi}$ the corresponding estimates.
The number of sampled units in the $lw$ group is denoted $m_{lw}$.
 Then, if the known proportion of units with $W=0$
 in the population is (say) 0.5, we may add 
 to our quadratic programming problem the constraint:
\begin{equation}
\sum_{l=1}^L  m_{l0} \sum_i \frac{{g}_{l0i}}{s_i} = 0.5I,
\end{equation}
Similarly for any explanatory variable with known population's proportion.
Such additional linear constraints  may resolve problems due
to non-identifiability, discussed in  subsection 2.5.

\bigskip

\subsection{ Approximation of the MSE of the modifeied Horvitz-Thompson
Estimator.}

In this subsection we will briefly discuss the estimation of the MSE of
our estimator $\hat{T}_l$, as defined in ($\ref{eqn:TL}$). 
The population total,
of items with a corresponding $X=l$, is denoted $T_l$.
The MSE of $T_1$
equals:

\begin{eqnarray}
E(\hat{T}_1-T_1)^2&=& E( \; (\hat{T}_1-\frac{N}{I}\theta_1)+ (\frac{N}{I}\theta_1 - T_1) \;)^2\\
&=&E (\hat{T}_1-\frac{N}{I}\theta_1)^2+ E(\frac{N}{I}\theta_1 - T_1) ^2.
\end{eqnarray}

In the above we represent the MSE as the sum of two terms,
the second term is  the variance of the (non-modified) Horvitz-Thompson estimator
$\frac{N}{I}\theta_1$. 

Note, our modifeid Horvitz-Thompson
procedure is essential when the non-response and whence
the bias is significant. In those cases the first term 
in the above equation is likely to dominate. This is also supported by our simulations, where it can be checked that the MSE
of the `non-modified' (or, oracle's version of) Horvitz-Thompson
is typically much smaller than that of the Modified Horvitz-Thompson.
 
The estimation, of the first term, may be obtained through
a parametric bootstrap, where we treat $\hat{G}_1$ as `true' and 
and we replace $T_1$ by $\frac{N}{I} E_{\hat{G}_1} \frac{1}{S}$.
We simulate $K$ times
$Y_1^k,...,Y_{m_1}^k$, $k=1,...,K$ according to $\hat{G}_1$,
then derive the
corresponding $\hat{T}^k_1$ and finally consider
$\sum_k \frac{(\hat{T}_1^k- \frac{N}{I} E_{\hat{G}_1} \frac{1}{S})^2}{K}$,
to obtain  an estimate of the first term.

It seems that our ideas may be applied also for the estimation
of the variance of the non-modified/oracle's Horvitz-Thompson estimator, at least in special cases, e.g.,  where the initial list is obtained through a simple random sample. However, since the estimation of
the variance of Horvitz-Thompson may be problematic even when
the sampling probabilities are exactly known, and in light of the
apparent dominance of the first term, we do not elaborate on it.

\subsection{ On the approximation of
$E_{{G}} \frac{1}{S}$ by $E_{\hat{G}} \frac{1}{S}$. }

From Lindsay and Roeder (1993), it follows  that 
in our truncated geometric example, where $Y$ may have only
$M_0$ distinct values, 
the mle estimate for $G$ 
is supported
on $M_0$ points at most. This implies that $\hat{G}$ may converge 
weakly to $G$,
as the number of observations increases, only if we assume that $G$ is supported on
at most $M_0$ points. 

The weak convergence argument has a special appeal in problems where 
the size of the 
support of $Y$ increases with the number of observations.
Such problems are described in subsection 2.6, and also briefly
in the following.
In principle our deconvolution method may also be applied in cases where $Y$
is continuous,
by discretizing $Y$ through some bins. When the discretezation
gets finer as the number of observations is increased
the support of the discretized $Y$, increases and 
we may expect
that $\hat{G}$ will converge weakly to $G$.

However, we may still hope  for a good approximation of $E_G \frac{1}{S}$ by the above method,
also when the true $G$ has more than $M_0$ points. From an EB
perspective we see that in fact it is enough to approximate
the Bayes decisions $E_G( \frac{1}{S}| Y=y)$ for the observed values
of $Y$. Our procedure may also be viewed 
from an EB perspective,
as estimation of
$\sum_i E_G( \frac{1}{S}| Y=y)$ by $\sum_i E_{\hat{G}} ( \frac{1}{S}| Y=y_i)$.
Estimation of $E_G( \frac{1}{S}| Y=y),$ for the observed values $y$, is easier than that of estimating $G$. Moreover, there could be different $G$ with the same Bayesian decisions for the observed values of $Y$. Consider for example the familiar Bayes estimator 
$E_G(\lambda|Y=y)$ for $\lambda$ in the Poisson case,  
it depends only on the ratio $P_G(Y=y+1)/P_G(Y=y)$,
the last ratio may be the same on the observed values for different $G$. Thus, we may expect similar performance of two different
mle estimates of $G$. 

The non-identifiability of $G$ may also be resolved
using additional calibration constraints as explained in  subsection 2.3.

\subsection{ Further applications of deconvolution.}

Though our motivating example is of estimating $\theta$, our approach may be useful  in a general setup where there is an interest in estimating
$\sum_i  h(\eta_i)$, for unknown parameters $\eta_1,...,\eta_m$, and  a given bounded continuous function $h$; the
setup is where the estimation is  based on independent
observations $Y_i \sim F_{\eta_i}$, $i=1,...,m$,
and $m$ is large;
here $Y_i$ may be truncated observations. 

In principle our method may also be applied in cases where $Y$
is continuous,
by discretizing $Y$ through some bins.

Note, in examples where there is a pointwise unbiased estimator for $h(\eta_i)$, $i=1,...m$,
plugging-in that pointwise unbiased estimator could be a preferred simple alternative to our approach of estimating first the empirical distribution $G$ of $\eta_1,...,\eta_m$. 
For example
consider the case where $F_{\eta_i}= N(\eta_i,1)$  and $h(\eta_i)=\eta_i$, $i=1,...,m$ (the observations are not truncated).

A slight twist of the last normal example with $F_{\eta_i}=N(\eta_i,1)$ 
is when truncation is involved,
e.g., when we want to estimate $\sum_{ \{ i: Y_i>T \} } \eta_i $. Here, our deconvolution
method is useful; we estimate the empirical distribution of 
the parameters $\eta_j$,
$ j \in \{ i: Y_i>T \} $, treating the remaining observations as truncated. Estimating the total amount of signal that corresponds to high valued observations might be of a special interest in multiple testing where we treat, the subset
$ \{ i: Y_i>T \} $, as potential discoveries. See, the treatment in Greenshtein et. al. (2008).

A related example, under the same setup,
is where $h(\eta_i)$ is an indicator of $\eta_i>C$, namely we want to estimate 
$\sum_{ \{ i: Y_i>T, \eta_i>C \} } h(\eta_i)=\# \{ i: Y_i>T, \eta_i>C \} $. The motivation for the last  estimation problem is related to
estimation of the False Discovery Rate, see Benjamini and Hochberg (1995). Here we think of $\eta_i>C$, for appropriate $C > 0$, as a meaningful discovery (rather than a `significant' discovery!), while 
$T$ is a threshold value for rejecting the null: $\eta_i \leq 0$. Among all the rejections/potential-discoveries, we want to estimate the amount of `true'/`meaningful' discoveries.
Again, we treat the non-rejected observations, i.e., those that correspond to indices $i$ with $Y_i \leq T$, as truncated. One may adjust $T$ 
(or in general the rejection region) to obtain a desired estimated false discovery rate. Our approach
is appropriate only asymptotically when $m$ is large.

A treatment of the problem of estimating random sums in other related examples, including a study of the efficiency of the
estimators, may be found in Zhang (2005).

\section{Simulations}

In the following reported simulations, we estimate a given proportion $pr_1$, based on
a random sample, as described above, with a list of size $I$.
We present simulations for three pairs of distributions of response probabilities. Each pair is consisted of $G_1$, the distribution of $\tilde{p}$ for items with $X=1$, and $G_0$, the distribution of 
$\tilde{p}$ for items with $X=0$. Given the list size $I$ we simulate
$I^1$ items in the list with corresponding $X_i=1$ and $I^0$ with corresponding $X_i=0$, $I^1 \sim B(I,pr_1)$, $I^0=I-I^1$.
For each item $i$ in the `$I^1$-list' we simulate a response probability $\tilde{p}_{1i}$ according to the distribution $G_1$; similarly for the other list we simulate response probabilities according to $G_0$. Finally, we simulate a 
(truncated by $M_0$) Geometric r.v., with the corresponding
$\tilde{p}_{li}$ for every item $i$ in the list $I^l$,
$l=0,1$. We stress again, in this section we describe
the simulations  in terms of $\tilde{p}_{li}$, and not in terms of $p^*_{li}$. 

The 3 pairs of distributions $(G_0,G_1)$ that we consider, are the following.

{\bf Two Points.}  The distribution $G_0$  has a two points support
at the points $0.5 $ and $0.9$, with probability  mass 0.5 at each.

The distribution $G_1 \equiv G_1^\alpha $ is a $(-\alpha)$ translation of $G_0$.
We present results for the cases $\alpha= 0.1, 0.2, 0.3, 0.4$.

\bigskip

{\bf Uniform.} The distribution $G_0$ is uniform on the interval $(0.1,1)$.
The distribution $G_1 \equiv G_1^\alpha$, is a mixture of $G_0$ and a point mass at $0.1$, where the mixing weights are $(1-\alpha)$ and $\alpha$ correspondingly. We present results for $\alpha=0.1, 0.2, 0.3, 0.4$. 

\bigskip

{\bf Normal} The distribution $G_0$ is a $N(0.5,0.1)$, left censored by
$0.1$ and right censored by $1$. The distribution $G_1 \equiv G_1^\alpha$ is 
$N(0.5-\alpha, 0.1)$ truncated bellow by $0.1$ and truncated above by $1$.
We present results for $\alpha=0.1, 0.2, 0.3, 0.4$. 

\bigskip

In all the above examples, as $\alpha$ increases the difference between $G_0$
and $G_1$ gets larger and
the treatment of the non-response observations as missing at random is less appropriate. Each pair of distributions was simulated for
the truncations $M_0=4,5,6,7$.

In all of our simulations we took a grid $\hat{\bdf{s}}$ whose points $\hat{s}_i$
correspond to
the choice of the $\tilde{p}_{l,i}$ points: $0.1, 0.12, 0.14,...,1$. The results are mildly sensitive to the choice of the grid. We used the method of moments version  coupled with the estimator (\ref{eqn:MHT1}) in our estimation,
the method will be reffered to as the Modified Horvitz Thompson
(MHT).

The following Table 1 and Table 2, summarize a
simulation study of the estimator (\ref{eqn:MHT1}) for the cases $I=1000$ and $I=10000$.
We kept $pr_1=0.5$.
Each entry is based on a simulation of 1000 repetitions. We present in the tables a comparison with the following two estimators. The naive-estimator,
that treats non-response as missing at random and estimate $pr_1$ by $m_1/(m_1+m_0)$; the oracle-estimator, see equation (\ref{eqn:HT}), which corresponds to the
Horvitz-Thompson estimator that an oracle, who knows the actual values of $p_i$,
would have used.
The  M-NV and M-MHT columns are the means of the naive estimator and the MHT estimator. One may see the bias correction of our MHT estimator relative to the naive. The oracle-estimator is obviously unbiased.
The S-NV, S-MHT and S-OR columns are the squared-roots of the MSE of the naive, Modified Horvitz-Thompson and the oracle estimators, correspondingly. The columns $\text{M-}m_0$ and $\text{M-}m_1$ give the mean
of $m_0$ and $m_1$, in the 1000 simulations.
As $\alpha$ increases, the advantage of our modification is more noticeable.
The variance of the naive estimator is smaller than that of the other two,
hence  in case $I=1000$ coupled with the small value $\alpha=0.1$, where the component of the  bias 
in the MSE is less dominant relative to that of the variance, 
the naive estimator occasionally performs better than both the modified Horvitz-Thompson and the Oracle estimators.
However, when $I=10000$ the oracle always dominates the naive.
In both cases, $I=1000$ and $I=10000$, our modified  Horvitz-Thompson 
estimator dominates the naive when $\alpha$ is large enough
and thus, it pays to reduce the bias on the expense of increasing the variance. As expected the oracle estimator always dominates our MHT estimator.  

The bias correction of our modified
Horvitz-Thompson procedure has merits beyond the reduction of the
squared loss, even if occasionally
the bias reduction is on the expense of increasing the overall squared risk. Often, the estimators arrive in a time series, e.g.,
proportion of unemployed in month $t$, and the final estimators involve smoothing of the time series. Then, smoothing around the
true value gives further squared loss reduction, compared to smoothing of a biased  sequence of estimators. 

\bigskip

{\it Note:} It seems that increasing $M_0$ is more effective in reducing the squared loss relative to  increasing  in $I$, when accounting for the extra effort in 
visits used by each strategy.


\begin{table}[htbp]
  \centering
  \caption{\bf{Simulations, I=1000}}
    \begin{tabular}{rrrrrrrrrr}
 $ \bdf{G_0}$&$\bdf{M_0}$ & $ \bdf{\alpha}$ & \textbf{M-NV} & \textbf{M-MHT} & \textbf{S-NV} & \textbf{S-MHT} & \textbf{S-OR} & $\bdf{\textbf{M-}m_1}$ & $\bdf{\textbf{M-}m_0}$ \\
    2points & 4     & 0.1   & 0.4953 & 0.5107 & 0.0059 & 0.0429 & 0.0076 & 472   & 481 \\
    2points & 4     & 0.2   & 0.4556 & 0.4823 & 0.0446 & 0.0508 & 0.0214 & 421   & 503 \\
    2points & 4     & 0.3   & 0.4515 & 0.4946 & 0.0488 & 0.0540 & 0.0142 & 395   & 480 \\
    2points & 4     & 0.4   & 0.4009 & 0.4556 & 0.0993 & 0.0707 & 0.0216 & 322   & 482 \\
    2points & 5     & 0.1   & 0.4760 & 0.4862 & 0.0242 & 0.0282 & 0.0186 & 463   & 510 \\
    2points & 5     & 0.2   & 0.4696 & 0.4949 & 0.0306 & 0.0311 & 0.0140 & 446   & 504 \\
    2points & 5     & 0.3   & 0.4697 & 0.5105 & 0.0306 & 0.0380 & 0.0166 & 426   & 480 \\
    2points & 5     & 0.4   & 0.4136 & 0.4637 & 0.0866 & 0.0516 & 0.0196 & 346   & 490 \\
    2points & 6     & 0.1   & 0.5200 & 0.5298 & 0.0201 & 0.0337 & 0.0242 & 512   & 472 \\
    2points & 6     & 0.2   & 0.4908 & 0.5139 & 0.0097 & 0.0253 & 0.0072 & 474   & 492 \\
    2points & 6     & 0.3   & 0.4553 & 0.4955 & 0.0449 & 0.0269 & 0.0144 & 424   & 507 \\
    2points & 6     & 0.4   & 0.4255 & 0.4750 & 0.0747 & 0.0377 & 0.0170 & 365   & 493 \\
    2points & 7     & 0.1   & 0.5125 & 0.5197 & 0.0126 & 0.0222 & 0.0152 & 508   & 483 \\
    2points & 7     & 0.2   & 0.4825 & 0.5010 & 0.0176 & 0.0156 & 0.0091 & 472   & 506 \\
    2points & 7     & 0.3   & 0.4651 & 0.5024 & 0.0351 & 0.0206 & 0.0114 & 440   & 506 \\
    2points & 7     & 0.4   & 0.4161 & 0.4623 & 0.0840 & 0.0432 & 0.0217 & 366   & 514 \\
    unif  & 4     & 0.1   & 0.4839 & 0.4941 & 0.0232 & 0.0602 & 0.0206 & 408   & 435 \\
    unif  & 4     & 0.2   & 0.4681 & 0.4864 & 0.0365 & 0.0650 & 0.0239 & 382   & 434 \\
    unif  & 4     & 0.3   & 0.4495 & 0.4786 & 0.0536 & 0.0678 & 0.0249 & 355   & 435 \\
    unif  & 4     & 0.4   & 0.4311 & 0.4753 & 0.0713 & 0.0702 & 0.0269 & 329   & 435 \\
    unif  & 5     & 0.1   & 0.4861 & 0.4923 & 0.0217 & 0.0424 & 0.0193 & 426   & 451 \\
    unif  & 5     & 0.2   & 0.4704 & 0.4901 & 0.0343 & 0.0451 & 0.0217 & 401   & 451 \\
    unif  & 5     & 0.3   & 0.4557 & 0.4839 & 0.0475 & 0.0482 & 0.0230 & 377   & 451 \\
    unif  & 5     & 0.4   & 0.4388 & 0.4793 & 0.0636 & 0.0533 & 0.0242 & 352   & 451 \\
    unif  & 6     & 0.1   & 0.4872 & 0.4961 & 0.0213 & 0.0334 & 0.0192 & 439   & 462 \\
    unif  & 6     & 0.2   & 0.4748 & 0.4914 & 0.0301 & 0.0362 & 0.0196 & 417   & 461 \\
    unif  & 6     & 0.3   & 0.4600 & 0.4871 & 0.0436 & 0.0402 & 0.0217 & 394   & 462 \\
    unif  & 6     & 0.4   & 0.4453 & 0.4843 & 0.0575 & 0.0413 & 0.0233 & 371   & 462 \\
    unif  & 7     & 0.1   & 0.4886 & 0.4955 & 0.0197 & 0.0279 & 0.0178 & 449   & 470 \\
    unif  & 7     & 0.2   & 0.4765 & 0.4926 & 0.0292 & 0.0311 & 0.0201 & 428   & 470 \\
    unif  & 7     & 0.3   & 0.4638 & 0.4894 & 0.0399 & 0.0319 & 0.0201 & 407   & 471 \\
    unif  & 7     & 0.4   & 0.4509 & 0.4892 & 0.0519 & 0.0328 & 0.0214 & 386   & 470 \\
    norm  & 4     & 0.1   & 0.4581 & 0.4865 & 0.0423 & 0.0533 & 0.0235 & 406   & 480 \\
    norm  & 4     & 0.2   & 0.4468 & 0.5107 & 0.0537 & 0.0644 & 0.0153 & 369   & 457 \\
    norm  & 4     & 0.3   & 0.3770 & 0.4992 & 0.1234 & 0.0694 & 0.0234 & 284   & 470 \\
    norm  & 4     & 0.4   & 0.3145 & 0.4866 & 0.1858 & 0.0611 & 0.0284 & 215   & 469 \\
    norm  & 5     & 0.1   & 0.4660 & 0.4833 & 0.0343 & 0.0351 & 0.0205 & 433   & 496 \\
    norm  & 5     & 0.2   & 0.4409 & 0.4896 & 0.0594 & 0.0414 & 0.0193 & 388   & 492 \\
    norm  & 5     & 0.3   & 0.3766 & 0.4734 & 0.1237 & 0.0526 & 0.0359 & 306   & 507 \\
    norm  & 5     & 0.4   & 0.3360 & 0.4855 & 0.1643 & 0.0461 & 0.0258 & 248   & 490 \\
    norm  & 6     & 0.1   & 0.4673 & 0.4802 & 0.0329 & 0.0275 & 0.0226 & 446   & 508 \\
    norm  & 6     & 0.2   & 0.4554 & 0.4923 & 0.0448 & 0.0280 & 0.0141 & 415   & 497 \\
    norm  & 6     & 0.3   & 0.4279 & 0.5056 & 0.0724 & 0.0355 & 0.0162 & 360   & 481 \\
    norm  & 6     & 0.4   & 0.3769 & 0.5079 & 0.1235 & 0.0376 & 0.0226 & 290   & 479 \\
    norm  & 7     & 0.1   & 0.5211 & 0.5304 & 0.0213 & 0.0333 & 0.0296 & 505   & 464 \\
    norm  & 7     & 0.2   & 0.4590 & 0.4872 & 0.0413 & 0.0244 & 0.0163 & 429   & 506 \\
    norm  & 7     & 0.3   & 0.4411 & 0.5048 & 0.0593 & 0.0289 & 0.0145 & 384   & 486 \\
    norm  & 7     & 0.4   & 0.4147 & 0.5282 & 0.0857 & 0.0411 & 0.0342 & 329   & 464 \\
    \end{tabular}%
  \label{tab:addlabel}%
\end{table}%


\begin{table}[htbp]
  \centering
  \caption{\bf{Simulations, I=10000}}
    \begin{tabular}{rrrrrrrrrr}
    $ \bdf{G_0}$&$\bdf{M_0}$ & $ \bdf{\alpha}$ & \textbf{M-NV} & \textbf{M-MHT} & \textbf{S-NV} & \textbf{S-MHT} & \textbf{S-OR} & $\bdf{\textbf{M-}m_1}$ & $\bdf{\textbf{M-}m_0}$ \\
    2points & 4     & 0.1   & 0.4928 & 0.5135 & 0.0073 & 0.0288 & 0.0026 & 4689  & 4827 \\
    2points & 4     & 0.2   & 0.4792 & 0.5221 & 0.0208 & 0.0397 & 0.0052 & 4418  & 4801 \\
    2points & 4     & 0.3   & 0.4499 & 0.5080 & 0.0501 & 0.0374 & 0.0054 & 3936  & 4813 \\
    2points & 4     & 0.4   & 0.4050 & 0.4671 & 0.0950 & 0.0465 & 0.0099 & 3249  & 4774 \\
    2points & 5     & 0.1   & 0.5011 & 0.5140 & 0.0014 & 0.0202 & 0.0073 & 4873  & 4852 \\
    2points & 5     & 0.2   & 0.4824 & 0.5171 & 0.0176 & 0.0271 & 0.0025 & 4580  & 4914 \\
    2points & 5     & 0.3   & 0.4622 & 0.5194 & 0.0378 & 0.0311 & 0.0058 & 4193  & 4879 \\
    2points & 5     & 0.4   & 0.4107 & 0.4750 & 0.0893 & 0.0341 & 0.0062 & 3438  & 4933 \\
    2points & 6     & 0.1   & 0.5024 & 0.5097 & 0.0025 & 0.0123 & 0.0065 & 4946  & 4898 \\
    2points & 6     & 0.2   & 0.4882 & 0.5133 & 0.0118 & 0.0196 & 0.0025 & 4719  & 4946 \\
    2points & 6     & 0.3   & 0.4639 & 0.5133 & 0.0361 & 0.0221 & 0.0039 & 4314  & 4985 \\
    2points & 6     & 0.4   & 0.4251 & 0.4873 & 0.0750 & 0.0204 & 0.0057 & 3649  & 4936 \\
    2points & 7     & 0.1   & 0.4976 & 0.5019 & 0.0025 & 0.0045 & 0.0008 & 4931  & 4979 \\
    2points & 7     & 0.2   & 0.4883 & 0.5053 & 0.0118 & 0.0106 & 0.0026 & 4773  & 5002 \\
    2points & 7     & 0.3   & 0.4692 & 0.5096 & 0.0308 & 0.0166 & 0.0046 & 4437  & 5019 \\
    2points & 7     & 0.4   & 0.4341 & 0.4926 & 0.0659 & 0.0135 & 0.0054 & 3802  & 4957 \\
    unif  & 4     & 0.1   & 0.4847 & 0.4938 & 0.0163 & 0.0399 & 0.0066 & 4084  & 4342 \\
    unif  & 4     & 0.2   & 0.4677 & 0.4894 & 0.0328 & 0.0444 & 0.0071 & 3818  & 4345 \\
    unif  & 4     & 0.3   & 0.4502 & 0.4870 & 0.0501 & 0.0452 & 0.0079 & 3557  & 4344 \\
    unif  & 4     & 0.4   & 0.4313 & 0.4831 & 0.0689 & 0.0473 & 0.0083 & 3293  & 4343 \\
    unif  & 5     & 0.1   & 0.4860 & 0.4944 & 0.0150 & 0.0285 & 0.0063 & 4262  & 4508 \\
    unif  & 5     & 0.2   & 0.4715 & 0.4932 & 0.0290 & 0.0304 & 0.0069 & 4019  & 4505 \\
    unif  & 5     & 0.3   & 0.4553 & 0.4888 & 0.0450 & 0.0338 & 0.0072 & 3768  & 4508 \\
    unif  & 5     & 0.4   & 0.4389 & 0.4855 & 0.0613 & 0.0351 & 0.0074 & 3524  & 4505 \\
    unif  & 6     & 0.1   & 0.4875 & 0.4964 & 0.0136 & 0.0203 & 0.0061 & 4394  & 4619 \\
    unif  & 6     & 0.2   & 0.4741 & 0.4930 & 0.0265 & 0.0224 & 0.0064 & 4165  & 4621 \\
    unif  & 6     & 0.3   & 0.4603 & 0.4928 & 0.0400 & 0.0243 & 0.0066 & 3939  & 4618 \\
    unif  & 6     & 0.4   & 0.4455 & 0.4885 & 0.0547 & 0.0271 & 0.0071 & 3711  & 4618 \\
    unif  & 7     & 0.1   & 0.4887 & 0.4967 & 0.0125 & 0.0158 & 0.0060 & 4493  & 4700 \\
    unif  & 7     & 0.2   & 0.4768 & 0.4953 & 0.0238 & 0.0177 & 0.0059 & 4282  & 4700 \\
    unif  & 7     & 0.3   & 0.4641 & 0.4929 & 0.0363 & 0.0194 & 0.0063 & 4072  & 4702 \\
    unif  & 7     & 0.4   & 0.4512 & 0.4919 & 0.0491 & 0.0196 & 0.0066 & 3865  & 4702 \\
    norm  & 4     & 0.1   & 0.4891 & 0.5110 & 0.0111 & 0.0321 & 0.0104 & 4327  & 4520 \\
    norm  & 4     & 0.2   & 0.4434 & 0.4993 & 0.0566 & 0.0373 & 0.0047 & 3670  & 4606 \\
    norm  & 4     & 0.3   & 0.3896 & 0.4973 & 0.1104 & 0.0433 & 0.0077 & 2922  & 4578 \\
    norm  & 4     & 0.4   & 0.3255 & 0.4944 & 0.1746 & 0.0376 & 0.0094 & 2209  & 4578 \\
    norm  & 5     & 0.1   & 0.4924 & 0.5079 & 0.0077 & 0.0196 & 0.0077 & 4565  & 4705 \\
    norm  & 5     & 0.2   & 0.4660 & 0.5090 & 0.0341 & 0.0257 & 0.0111 & 4082  & 4677 \\
    norm  & 5     & 0.3   & 0.4029 & 0.4888 & 0.0972 & 0.0320 & 0.0066 & 3245  & 4809 \\
    norm  & 5     & 0.4   & 0.3377 & 0.4787 & 0.1623 & 0.0347 & 0.0136 & 2491  & 4885 \\
    norm  & 6     & 0.1   & 0.4847 & 0.4963 & 0.0153 & 0.0124 & 0.0049 & 4621  & 4912 \\
    norm  & 6     & 0.2   & 0.4694 & 0.5023 & 0.0306 & 0.0176 & 0.0051 & 4273  & 4830 \\
    norm  & 6     & 0.3   & 0.4213 & 0.4908 & 0.0787 & 0.0232 & 0.0049 & 3553  & 4881 \\
    norm  & 6     & 0.4   & 0.3770 & 0.4988 & 0.1231 & 0.0217 & 0.0103 & 2899  & 4791 \\
    norm  & 7     & 0.1   & 0.4944 & 0.5028 & 0.0056 & 0.0084 & 0.0029 & 4790  & 4898 \\
    norm  & 7     & 0.2   & 0.4799 & 0.5067 & 0.0201 & 0.0141 & 0.0077 & 4476  & 4851 \\
    norm  & 7     & 0.3   & 0.4406 & 0.4983 & 0.0595 & 0.0172 & 0.0069 & 3834  & 4869 \\
    norm  & 7     & 0.4   & 0.3950 & 0.4993 & 0.1050 & 0.0177 & 0.0097 & 3164  & 4845 \\
    \end{tabular}%
  \label{tab:addlabel}%
\end{table}%

\newpage

\section{ Analysis of real data of Labor Force Survey.}  
  
In this section we will apply our method on real data from the
Labor Force Survey that is conducted by the Israel Central Bureau of Statistics.  The sampling method is 4-8-4 rotating panels, however
for our analysis, it may be equivalently trated and described as
a 4-in rotation, which is described in the following.

The survey is given to four panels, where each panel is investigated
for four consecutive months. Each month one panel finishes its fourth
investigation and in the next month it will be replaced by a new panel that will remain for four months. The main purpose of the survey is to estimate the proportion of  `Unemployment', `Employment',
and those who are `Not in Working Force (NWF)', the last category is of those that do not have a job nor they are looking for a job;
denote the corresponding values of our variable $X$-`working status', by
0, 1, 2.  We are interested in estimating $pr_0, pr_1$ and $pr_2$.
The population 
of interest is  of residents whose  age is above 15, and the proportions are with respect to that population.

Temporarily assume that, we have only the data from the panel that is investigated for the fourth
time. Its size is $I$ ( about 5000; the size of the entire list
of the four panels is about 20000 each month);  however, only $m$ responses were
obtained $m_l$ responses from people with working status $l$, $l=0,1,2$. 
The response rate is about 80 percent each month.
For each of those $m$ units there is a corresponding
truncated random variable, denoted $Y$, that counts the number of responses, 
including the current one. 
We model the distribution of this truncated random variable by
$$ Y=1+W; \;\; W \sim Binomial(3,p^*).$$

The above model amounts to assuming that the probability of response
of unit $i$,
is $p^*_i$ in all of its four investigation attempts, and responses in different months are independent. 
Given a grid
$\hat{s}_1,...,\hat{s}_{\kappa}$, a matrix $P_{\bf{\hat{s}}}=(p_{ji})$
is defined where $p_{ji}=  P_{\hat{s}_i}(Y=j)=P_{\hat{s}_i}(1+W=j)$,
$j=1,2,3,4$,
for $W \sim B(3,\hat{s}_i)$. In our analysis we took the grid
0.1, 0.11, 0.12,...,1. Thus, our matrix $P_{{\bf \hat{s}}}$ is 
of dimension $4 \times 90$.

Now, $pr_l$, $l=0,1,2$, may be estimated as in (\ref{eqn:MHT1}).
However, so far we used only the data from the panel that has four investigations. Indeed the panels that have less investigations
will yield poor estimates of $E_{G_l} 1/S$. Our approach is the following hybrid method in which we estimate $E_{G_l} 1/S$,
$l=0,1,2$, based on
the data from the current `fourth panel' in addition to 
the data obtained in the fourth investigation of the three more panels that had their fourth investigation in the previous month,
two months ago and three months ago, altogether four panels. Let
$m_0,m_1,m_2$, the number of items in the currently investigated four panels, with corresponding $X=l$, $l=0,1,2$.
Our hybrid approach is to inflate $m_l$ which is based on
the currently investigated four panels, using the estimated $E_{G_l} 1/S$,
$l=0,1,2$ which are based on current as well as `historical' complementary information.
The underlying assumption is that $E_{G_l} 1/S$  changes slowly in time and thus estimating it based on a complementary older data
we still  get at least some bias correction.

We finally get the  estimator  $$ \hat{pr}_l=\frac{ m_l E_{\hat{G}_l} \frac{1}{S}}{ \sum_l  m_l E_{\hat{G}_l} \frac{1}{S}}.$$

\bigskip

Since the true proportions of the various working stauses
are unknown, we will first demonstrate
the performance of the above estimation method in estimating the following {\it known} true proportions, based on the responses in a given month.

In one case  we  estimate the proportion of males in the population,
which is known to be 0.4853; their proportion in the survey
among responses
is about one percent lower. In the other example we estimate the proportion of the group age 20-39. Their known proportion is 0.397
while their, 
response rate is particularly low, their proportion  among the responses is nearly 3 percent lower than their proportion in the population. 

Each of the following tables 3 and 4 has three lines that correspond
to the  data obtained in Aug/2012, Dec/2012, and April/2013. We took
periods that are four months apart in order not to have overlapping
panels. The general picture persist in other  months.

The columns  True, Naive, and MHT correspond to the true proportion,
the sample proportion among responses, and our modified Horvitz Thompson estimator. In each case one may see that the MHT corrects
in the right direction.

\bigskip

After gaining some confidence in the MHT, we will now examine
its estimates in the estimation of the proportion of `Unemployed',
`Employed' and those `Not in Working Force' (NWF).
In the following Table 5 the columns Naive and MHT
are as before. The column Bureau gives the estimates of the 
Israel, Central Bureau of Statistics for the three categories of working status. The three parts of the table refer to the three working statuses. The three lines in each part refer to the three months as described before. The Bureau and the MHT  `correct'
the naive estimator for Employment and NWF,
in opposite directions  (the official Bureau estimator
involves additional seasonal adjustment that we neglect). The estimator of the bureau is obtained through a method that involves calibration
in a `post-stratification manner'. It seems that the correction of the bureau of  `Employment' and the `NWF' is in the wrong direction. This is indicated also when imputing missing values based
on their values in months where a response was obtained looking also `into the future'. On the other hand both the Bureau and the MHT estimators correct the unemployment naive estimate by increasing it.
This direction of correction of unemployment,
is also supported by an analysis that involves imputation.   

\newpage

\begin{table}[ht]
\caption{Comparison of estimates of male's proportion.} % title of Table
\centering % used for centering table
\begin{tabular}{c c c c} % centered columns (4 columns)
\hline\hline %inserts double horizontal lines
   & True & Naive & MHT \\ [0.5ex] % inserts table
%heading
\hline % inserts single horizontal line
{\text Male}&0.4853   & 0.4752 & 0.4822\\
            & 0.4853 & 0.4751 & 0.4819 \\ % inserting body of the table
            &0.4853 &  0.4776 & 0.4842 \\
%{\text UnEmp}& 0.3465 & 0.3576 & 0.3605 \\
% & 31 & 25 & 415 \\
%4 & 35 & 144 & 2356 \\
%{\text NWF} & 0.0454 & 0.0431 & 0.04840\\ [1ex] % [1ex] adds vertical space
\hline %inserts single line
\end{tabular}
\label{table:nonlin2} % is used to refer this table in the text
\end{table}

\begin{table}[ht]
\caption{ Comparison of estimates of  proportion of 20-39
age group.} % title of Table
\centering % used for centering table
\begin{tabular}{c c c c} % centered columns (4 columns)
\hline\hline %inserts double horizontal lines
   & True & Naive & MHT \\ [0.5ex] % inserts table
%heading
\hline % inserts single horizontal line
{\text Age 20-39} & 0.3970 & 0.3664 & 0.3815 \\
                  & 0.3970 & 0.3631 & 0.3984 \\ % inserting body of the table
                 & 0.3970 & 0.3598 &  0.3842\\
%{\text UnEmp}& 0.3465 & 0.3576 & 0.3605 \\
% & 31 & 25 & 415 \\
%4 & 35 & 144 & 2356 \\
%{\text NWF} & 0.0454 & 0.0431 & 0.04840\\ [1ex] % [1ex] adds vertical space
\hline %inserts single line
\end{tabular}
\label{table:nonlin3} % is used to refer this table in the text
\end{table}

\begin{table}[ht]
\caption{Comparison of unemployment estimates.} % title of Table
\centering % used for centering table
\begin{tabular}{c c c c} % centered columns (4 columns)
\hline\hline %inserts double horizontal lines
   & Bureau & Naive & MHT \\ [0.5ex] % inserts table
%heading
\hline % inserts single horizontal line
{\text Emp}
            & 0.6104 & 0.5931 & 0.5761\\
            & 0.6081 & 0.5992 & 0.5910 \\ % inserting body of the table
            &0.6089 & 0.5986 & 0.5881 \\
{\text NWF}&0.3416  & 0.3594 & 0.3748\\
           & 0.3465 & 0.3576 & 0.3605 \\
           &0.3491  & 0.3621 & 0.3720 \\
% & 31 & 25 & 415 \\
%4 & 35 & 144 & 2356 \\
{\text UnEmp} &0.0479  & 0.0475 & 0.0492 \\
              & 0.0454 & 0.0431 & 0.0484\\ 
              &0.0420 & 0.0392 & 0.0399 \\
              %[1ex] % [1ex] adds vertical space
               
\hline %inserts single line
\end{tabular}
\label{table:nonlin1} % is used to refer this table in the text
\end{table}

\newpage

\vspace{3ex}
{\bf \Large References}

\begin{list}{}{\setlength{\itemindent}{-1em}\setlength{\itemsep}{0.5em}}

\item

Benjamini, Y. and Hochberg, Y. (1995). Controling the false discovery rate: a practical and powerful approach to multiple testing. {\it JRSSB} {\bf 57} No.1, 289-300.

\item
Brown, L.D. and Greenshtein, E. (2009). Non parametric
empirical Bayes and compound decision
approaches to estimation of high dimensional vector of normal
means. {\it Ann. Stat.} {\bf 37}, No. 4, 1685-1704.

\item

Brown L.D., Greenshtein, E. and Ritov, Y. (2013). The Poisson compound decision revisited. {\it JASA.} {\bf 108} 741-749.

\item

Efron, B. (2003).  Robbins, Empirical Bayes and Microarrays.
{\it Ann.Stat.} {\bf 31} No. 2,  366-378.

%\item

%Efron B, and Thisted, R. (1976).  Estimating the number of unseen
%species: How many words did Shakespeare know? {\it Biometrika}
%{\bf 63} 435-447.

\item

Efron, B. (2013). Empirical Bayes modeling, computation and accuracy.
Manuscript.
\item

Greenshtein, E., Park, J., and Ritov, Y. (2008). Estimating the mean of high valued observations in high dimensions. {\it JSTP} {\bf 2} No. 3 407-418.

\item

Jiang, W. and Zhang, C.-H. (2009). General maximum likelihood
empirical Bayes estimation of normal means. {\it Ann. Stat.} {\bf 37}, No 4, 1647-1684.

\item
Koenker, R. and Mizera, I. (2013). Convex optimization, shape constraints, compound decisions and empirical Bayes rules. Manuscript.
\item

Lee, M., Hall, P., Haipeng, S., Marron, J.S., and Tolle, J. (2013).
Deconvolution estimation of mixture distributions with boundaries. (2013). {\it Electronic J. of Stat.} {\bf 7} 323-341.

\item

Lindsay, B. G. (1995). Mixture Models: Theory, Geometry and Applications.
Hayward, CA, IMS.

\item

Lindsay, B. G. and  Roeder, K. (1993).  Uniqueness of estimation and
identifiability in mixture models. {\it Canadian Journal of Stat.}
{\bf 21}, No. 2, 139-147.

%\item

%Link, W. A. (2004). Individual heterogeneity and identifiability in capture recapture models.  {\it Animal biodiversity and conservation} 27.1, 87-91.

\item
Little, R.J.A and Rubin, D.B.  (2002). Statistical Analysis with Missing Data. New York: Wiley.

\item

Karatzoglou, A., Smola, A., Hornik, K., and Zeleis, A., (2004).
An S4 package for kernel methods in R. {\it Journal of Statistical Software}
{\bf 11}, No. 9, 1-20.

\item

Kiefer and Wolfowitz (1956). Consistency of the maximum likelihood estimator in the presence  of infinitely many incidental parameters.
{\it Ann.Math.Stat.} {\bf 27} No. 4, 887-906.

%\item

%Norris, J. L. and Pollock, K. H. (1996). Nonparametric MLE under two closed capture-recapture
%models with heterogeneity. {\it Biometrics }, {\bf 52}, 639�649.

\item
Sharon L. Lohr (1999). \emph{Sampling Design and Analysis}. Brooks/Cole publishing company.

%\item

%Vapnik, N. V. (1998). Statistical Learning Theory. Wiley, New York.

\item

Zhang, C-H. (2005). Estimation of sums of random variables: Examples and information bounds. {\it Ann. Stat.} {\bf 33} No.5. 2022-2041.

\end{list}

\end{document}